\documentclass[12pt]{amsart}

\usepackage{amssymb}
\usepackage{latexsym}
\usepackage{amsmath}

\begin{document}


\def\c{{\mathbb C}}
\def\r{{\mathbb R}}
\def\q{{\mathbb Q}}
\def\z{{\mathbb Z}}
\def\i{{\mathbb I}}

\def\p{{\mathcal P}}
\def\h{{\mathcal H}}

\newtheorem{theo}{Theorem}[section]
\newtheorem{defi}{Definition}[section]
\newtheorem{coro}{Corollary}[section]
\newtheorem{exam}{Example}[section]
\newtheorem{rema}{Remark}[section]

\newtheorem{pro}{Problem}
\newtheorem{thm}{Theorem}
\newtheorem{cor}{Corollary}

\title{Quantum automorphism groups of homogeneous graphs}
\author{Teodor Banica}
\address{Institut de Mathematiques de Jussieu, 175 rue du Chevaleret, 75013 Paris}
\email{banica@math.jussieu.fr}
\maketitle


{\bf Abstract.} Associated to a finite graph $X$ is its quantum automorphism group $G$. The main problem is to compute the Poincar\'e series of $G$, meaning the series $f(z)=1+c_1z+c_2z^2+\ldots$ whose coefficients are multiplicities of $1$ into tensor powers of the fundamental representation. In this paper we find a duality between certain quantum groups and planar algebras, which leads to a planar algebra formulation of the problem. Together with some other results, this gives $f$ for all homogeneous graphs having $8$ vertices or less.


\section*{Introduction}

A remarkable discovery, due to Wang \cite{wa2}, is that the set $X_n=\{ 1,\ldots ,n\}$ has a quantum automorphism group, bigger in general than the symmetric group $S_n$. The quantum group doesn't exist of course, but the algebra of continuous functions on it does. This is a certain Hopf $\c^*$-algebra constructed with generators and relations, denoted here $H(X_n)$.

For $n=1,2,3$ the quotient map $H(X_n)\to\c (S_n)$ is an isomorphism. For $n\geq 4$ it is not, and in fact $H(X_n)$ is non commutative, and infinite dimensional.

There are several variations of this construction, see for instance Bichon \cite{bic} and \cite{sms}. The idea is that Hopf algebra quotients of $H(X_n)$ correspond to quantum automorphism groups of various discrete objects, like finite graphs, finite metric spaces, and so on.

A first purpose of this paper is to formulate some precise problem concerning such universal Hopf algebras. We believe that the following statement is the good one.

\begin{pro}
Let $X$ be a finite graph, all whose edges are colored and possibly oriented, such that an oriented edge and a non-oriented one cannot have same color.

We denote by $H(X)$ the universal Hopf $\c^*$-algebra coacting on $X$. This is obtained as an appropriate quotient of Wang's algebra $H(X_n)$, where $n$ is the number of vertices of $X$.

The problem is to compute its Poincar\'e series $f(z)=1+c_1z+c_2z^2+\ldots$, where $c_k$ is the multiplicity of $1$ into the $k$-th tensor power of the fundamental coaction of $H(X)$.
\end{pro}

As a first remark, in the computation $X\to f$ both the input and the output are classical objects. The problem is to find a classical computation relating them.

The input is a bit more general than in our previous paper \cite{sms}. Indeed, the algebra associated in \cite{sms} to a finite metric space $Y$ is the same as the algebra $H(Y_{gr})$, where $Y_{gr}$ is the complete graph having vertices at points of $Y$, with edges colored by corresponding lengths.

As for the output, the choice of $f$ instead of other invariants, like fusion rules etc. is inspired from recent progress in subfactors, such as Jones's fundamental work \cite{j3}.

There is actually one more subtlety here, concerning the input. The interesting case is when $X$ is quantum homogeneous, meaning that the algebra of fixed points of the universal coaction reduces to the scalars. (This is the same as asking for the equality $c_1=1$, known to correspond to the irreducibility condition in subfactors.) It follows from definitions that if $X$ is homogeneous then it is quantum homogeneous. The converse appears to be true in many cases of interest, but so far we don't know if it is true in general.

The main tool for solving problem 1 is Woronowicz's Tannakian duality \cite{wo2}. In this paper we find a general result in this sense. This is a duality between Hopf algebra quotients of $H(X_n)$ and subalgebras of Jones's spin planar algebra \cite{j1}, which preserves Poincar\'e series. The algebra $H(X)$ corresponds in this way to the planar algebra generated by the incidency matrices of $X$, one for each color, viewed as 2-boxes in the spin planar algebra.

The main application is with a product of $s$ complete graphs. The $s$ incidency matrices satisfy Landau's exchange relations \cite{la}, so we get a Fuss-Catalan algebra on $s$ colors, whose Poincar\'e series is computed by Bisch and Jones in \cite{bj1}. The $s=2$ case of this result was previously obtained in \cite{sms}, as a corollary of a direct combinatorial computation in \cite{fc}.

We get in this way a Fuss-Catalan series of graphs, containing products of complete graphs, plus other graphs, obtained by performing modifications which preserve the $H$ algebra.

Another series is the dihedral one, where we have $H(X)=\c (D_n)$. It is known from \cite{sms} that $n$-gons with $n\neq 4$ belong to it, and here we find an improved statement.

It is known from \cite{sms} that the Fuss-Catalan and dihedral series cover all non-oriented homogeneous graphs having $n\leq 7$ vertices. In this paper we prove the following result.

\begin{thm}
The non-colored non-oriented homogeneous graphs having $n\leq 8$ vertices fall into three classes: (1) Fuss-Catalan graphs, (2) dihedral graphs, (3) the cube and its complement, which correspond to a tensor product between $TL(2)$ and $TL(4)$.
\end{thm}

The case of colored or oriented graphs is discussed as well.

This work was done with help from several people, and in particular I would like to thank Dietmar Bisch, Ga\"etan Chenevier and Etienne Ghys.

\section{Formalism}

Let $H$ be a $\c^*$-algebra with unit, together with $\c^*$-morphisms $\Delta :H\to H\otimes H$,  $\varepsilon :H\to\c$ and $S:H\to H^{op}$ called comultiplication, counit and antipode. Here $\otimes$ is any $\c^*$-algebra tensor product and $H^{op}$ is the $\c^*$-algebra $H$, but with opposite product. We assume that the square of the antipode is the identity, and that Woronowicz's axioms in \cite{wo3} are satisfied.

Let $X$ be a finite set. We denote by $\c (X)$ the algebra of complex functions on $X$. The linear form on $\c (X)$ which sums the values of the function is denoted $\Sigma$.

\begin{defi}
A coaction of $H$ on $X$ is a morphism of $\c^*$-algebras $v :\c (X)\to \c (X)\otimes H$ satisfying the following conditions.

(a) Coassociativity condition $(id\otimes\Delta )v=(v\otimes id)v$.

(b) Counitality condition $(id\otimes\varepsilon )v=id$.

(c) Natural condition $(\Sigma\otimes id)v=\Sigma(.)1$.
\end{defi}

We should mention that in this definition the terminology is not standard.

The natural condition says that the action of the corresponding quantum group must preserve the counting measure on $X$. This condition is satisfied in all reasonable situations, but is not automatic in general. See Wang \cite{wa2}.

Consider the basis of $\c (X)$ formed by Dirac masses. Any linear map $v :\c (X)\to \c (X)\otimes H$ can be written in terms this basis, and the matrix of coefficients $(v_{ij})$ determines $v$.
$$v(\delta_i)=\sum_j \delta_j\otimes v_{ji}$$

The simplest example involves the Hopf $\c^*$-algebra $\c (G)$ associated to a finite group $G$. Here the comultiplication, counit and antipode are obtained by applying the $\c$ functor to the multiplication, unit and inverse map of $G$. If $G$ acts by permutations on $X$, we can apply the $\c$ functor to the corresponding map $(x,g)\to g(x)$ and we get a coaction of $\c (G)$ on $X$.

\begin{theo}
Let $G$ be a group of permutations of a finite set $X$. The coefficients $v_{ij}$ of the corresponding coaction of $\c (G)$ on $X$ are characteristic functions of the sets $\{ g\in G\mid g(j)=i\}$. When $i$ is fixed and $j$ varies, or vice versa, these sets form a partition of $G$.
\end{theo}

Back to the general case, it is convenient to translate the conditions in definition 1.1 in terms of coefficients $v_{ij}$. First, the Dirac masses being self-adjoint, the fact that $v$ is a $*$-map says that the elements $v_{ij}$ are self-adjoint. Multiplicativity of $v$ says that $v(\delta_i)v(\delta_k)=v(\delta_i\delta_k)$ for any $i,j,k$. This translates into the following formulae, where $\delta_{ik}$ is a Kronecker symbol.
$$\sum_j \delta_j\otimes v_{ji}v_{jk}
= \delta_{ik}\sum_j \delta_j\otimes v_{ji}$$

For $i=k$ we must have that $v_{ji}$ is a projection, and for $i\neq k$ we must have that $v_{ji}$ is orthogonal to $v_{jk}$. The fact that $v$ is unital translates into the following formula.
$$\sum_{ij} \delta_j\otimes v_{ji}=\sum_j\delta_j\otimes 1$$

Summing up, the fact that the linear map $v$ is a morphism of $\c^*$-algebras is equivalent to the fact that all rows of the matrix $v$ are partitions of unity with self-adjoint projections. The natural condition says that the sum on each column of $v$ is 1. Thus $v$ is a unitary corepresentation, so the antipode is given by $S(v_{ij})=v_{ji}$. By applying the antipode we get that the columns of $v$ are also partitions of unity. In other words, the linear map produced by a matrix $v$ is a coaction if and only if the coassociativity, counitality and natural conditions are satisfied, and $v$ is a ``magic biunitary'' in the following sense.

\begin{defi}
A matrix $v\in M_X(H)$ is called magic biunitary if its rows and columns are partitions of unity of $H$ with self-adjoint projections.
\end{defi}

A magic biunitary $v$ is indeed a biunitary, in the sense that both $v$ and its transpose $v^t$ are unitary matrices. The terminology comes from a vague similarity with magic squares.

A coaction is said to be faithful if its coefficients generate the $\c^*$-algebra $H$. We are interested in faithful coactions, and it is convenient to translate both axioms for $H$ and for $v$ in terms of coefficients, by formulating the above discussion in the following way.

\begin{theo}
If $v$ is a faithful coaction of $H$ on $X$ the following conditions are satisfied.

(a) The matrix $v=(v_{ij})$ is a magic biunitary and its coefficients generate $H$.

(b) There is a $\c^*$-morphism  $\Delta :H\to H\otimes H$ such that $\Delta (v_{ij})=\sum v_{ik}\otimes v_{kj}$.

(c) There is a $\c^*$-morphism $\varepsilon :H\to\c$ such that $\varepsilon (v_{ij})=\delta_{ij}$.

(d) There is a $\c^*$-morphism $S:H\to H^{op}$ such that $S(v_{ij})=v_{ji}$.

Conversely, a $\c^*$-algebra $H$ together with a matrix $v\in M_X(H)$ satisfying these conditions has a unique Hopf $\c^*$-algebra structure such that $v$ is a faithful coaction of $H$ on $X$.
\end{theo}

This statement is similar to definition 1.1 in Woronowicz's paper \cite{wo1}. The objects $(H,v)$ satisfying the above conditions correspond to ``compact permutation pseudogroups''.

Woronowicz's analogue of the Peter-Weyl theorem in \cite{wo1} shows that each irreducible corepresentation of $H$ appears in a tensor power of $v$. The main problem is to decompose these tensor powers, given by the following formulae.
$$v^{\otimes k}(\delta_{i_1}\otimes \ldots\otimes\delta_{i_k})
=\sum(\delta_{j_1}\otimes \ldots\otimes\delta_{j_k})
\otimes v_{j_1i_1}\ldots v_{j_ki_k}$$

By Frobenius reciprocity an equivalent problem is to compute spaces of fixed points.
$$Hom(1,v^{\otimes k})=\{ x\in\c (X)^{\otimes k}\mid v^{\otimes k}(x)=x\otimes 1\}$$

A slightly easier problem is to compute dimensions, arranged in a Poincar\'e type series.

\begin{defi}
The Poincar\'e series of a coaction $v:\c (X)\to\c (X)\otimes H$ is
$$f(z)=\sum_{k=0}^\infty dim\left(Hom(1,v^{\otimes k})\right)z^k$$
where $Hom(1,v^{\otimes k})$ is the space of fixed points of the $k$-th tensor power of $v$.
\end{defi}

Woronowicz's results in \cite{wo1} show that the dimension of the fixed point space of a corepresentation is obtained by applying the Haar functional $h:H\to\c$ to its character. By using the multiplicativity of the character map $\chi$ we get the following equality.
$$dim\left(Hom(1,v^{\otimes k})\right) =h\left(\chi_v^k\right)$$

The diagonal entries of $v$ are self-adjoint projections, so their sum $\chi_v$ is self-adjoint. For $z$ small enough the operator $1-z\chi_v$ is invertible, with inverse given by the following formula.
$$(1-z\chi_v)^{-1}=\sum_{k=0}^\infty z^k\chi_v^k$$

By the above, applying the Haar functional gives the Poincar\'e series.

\begin{theo}
The Poincar\'e series of a coaction $v:\c (X)\to\c (X)\otimes H$ is given by
$$f(z)=h\left( \left(1-z\sum_iv_{ii}\right)^{-1}\right)$$
where $h:H\to\c$ is the Haar functional. If $v$ is faithful the convergence radius of $f$ is bigger than $1/(\# X)$, with equality if and only if $H$ is amenable in the Hopf $\c^*$-algebra sense.
\end{theo}

In this statement the first assertion follows from the above discussion. The second one follows from a quantum analogue of the Kesten amenability criterion for discrete groups, due to Skandalis, and written down in the last section of \cite{subf}.  

For $H=\c (G)$ the Poincar\'e series is given by the following formula, where $\pi$ is the corresponding representation of $G$ on $\c (X)$, and where $\chi_\pi$ is its character.
$$f(z)=\frac{1}{|G|}\sum_{g\in G}\frac{1}{1-\chi_\pi (g)z}$$

The number $\chi_\pi (g)$ is the trace of $g$, viewed as an operator from $\c (X)$ to itself. The Hilbert space structure of $\c (X)$ is the one making the Dirac masses an orthonormal basis. Since $g$ permutes elements of the basis, its trace $\chi_\pi (g)$ is the number $m_g$ of its fixed points.

\begin{theo}
Let $G$ be a group of permutations of a set $X$ with $n$ elements. The Poincar\'e series of the corresponding coaction of $\c (G)$ is given by the formula
$$f(z)=\frac{1}{|G|}\sum_{m=0}^n\frac{\# G_m}{1-mz}$$
where $G_m\subset G$ is the subset of permutations having exactly $m$ fixed points. The convergence radius is $1/n$.
\end{theo}

Trying to compute Poincar\'e series will be the main objective in this paper. For a noncommutative Hopf $\c^*$-algebra $H$ this is in general a quite complicated analytic function, and doesn't have such a simple decomposition as a sum. For instance in the non-amenable case computation of the convergence radius is known to be a delicate problem.

\section{Coactions on abstract spaces}

Let $X$ be a finite set and let $d=(d_{ij})$ be a complex matrix with indices in $X$.  If a group $G$ acts on $X$ as in section 1 then both $d$ and $v=(v_{ij})$ are matrices with coefficients in $H$ and indices in $X$, so we can form the products $dv$ and $vd$.

\begin{theo}
Let $G$ be a group of permutations of a finite set $X$ and let $d\in M_X(\c)$ be a complex matrix with indices in $X$. Consider the magic biunitary $v$ describing the corresponding coaction of $\c (G)$. The action of $G$ preserves coefficients of $d$ if and only if $dv=vd$.
\end{theo}

\begin{proof}
The product $dv$ is computed by using the formula for $v_{ij}$ in section 1. 
$$(dv)_{ij}=\sum_kd_{ik}v_{kj}=\sum_k\sum_{g(j)=k}d_{ik}\delta_g=\sum_gd_{ig(j)}\delta_g$$

The same method gives a similar formula for the product $vd$. 
$$(vd)_{ij}=\sum_kv_{ik}d_{kj}=\sum_k\sum_{g(k)=i}d_{kj}\delta_g=\sum_gd_{g^{-1}(i)j}\delta_g$$

With $i=g(l)$ we get that $d_{lj}=d_{g(l)g(j)}$ holds for any $l,j,g$ if and only if $dv=vd$.
\end{proof}

\begin{exam}
Finite metric spaces.
\end{exam}

Let $X$ be a finite metric space. That is, we are given a finite set
$X$ and a real function $d:X\times X\to\r$ which is zero on the
diagonal, positive outside, and whose values satisfy the triangle
inequality. The distance function $d$ can be regarded as a complex
matrix with indices in $X$. Then the action of $G$ is isometric if and
only if $dv=vd$.

This follows from theorem 2.1. Indeed, the action is isometric
when $d_{ij}=d_{g(i)g(j)}$ for any $i,j\in X$, and this
means that $G$ preserves coefficients of $d$ in the sense of theorem 2.1.

\begin{exam}
Finite graphs.
\end{exam}

Let $X$ be a finite graph. That is, we are given a finite set $X$, whose elements are called vertices, and edges are drawn between certain pairs of different vertices. The edges are uniquely determined by the incidency matrix, given by $d_{ij}=1$ if $i,j$ are connected by an edge and $d_{ij}=0$ if not. Then the action of $G$ preserves the edges if and only if $dv=vd$.

This is another application of theorem 2.1. For, recall first that $G$
preserves the edges of a graph with vertex set $X$ when $(ij)$ is an edge $\Leftrightarrow$
$(g(i)g(j))$ is an edge. In terms of the incidency matrix, this
condition is $d_{ij}=1$ $\Leftrightarrow$
$d_{g(i)g(j)}=1$. Now since $d$ is a $0-1$ matrix this is the same as
asking for the equalities $d_{ij}=d_{g(i)g(j)}$ for any $i,j\in X$, so 
theorem 2.1 applies.

\begin{exam}
Finite oriented graphs.
\end{exam}

Let $X$ be a finite oriented graph. That is, $X$ is a graph all whose
edges have an orientation. We associate to $X$ matrices $r$ and $s$ as
follows. If $ij$ is an oriented edge we set $r_{ij}=1$, $r_{ji}=0$ and
$s_{ij}=i$, $s_{ji}=-i$. These matrices are related by the formula
$s=ir-ir^t$. Each of them determines $X$. The matrix $r$ has the
advantage of being real and the matrix $s$ has the advantage of being
self-adjoint. The conditions $rv=vr$ and $sv=vs$ are equivalent, and
are satisfied if and only if the action of $G$ preserves the oriented
edges.

This is seen as folows. First, the action of $G$ preserves the edges
of an oriented graph with vertex set $X$ when $(ij)$ is an oriented edge $\Leftrightarrow$
$(g(i)g(j))$ is an oriented edge.

In terms of $r$ this is equivalent
to $r_{ij}=1$ $\Leftrightarrow$ $r_{g(i)g(j)}=1$, and since $r$ is a $0-1$ matrix this is the same as
asking for the equalities $r_{ij}=r_{g(i)g(j)}$ for any $i,j\in
X$. Thus theorem 2.1 applies and shows that $G$ preserves the oriented
edges if and only if $rv=vr$.

In terms of $s$ we have that $G$ preserves the oriented edges if and
only if $s_{ij}=\pm i$ $\Leftrightarrow$
$s_{g(i)g(j)}=\pm i$, for any $i,j\in X$ and for any
choice of the sign $\pm$. Since $\pm i$ are the only non-zero
coefficients of $s$, this is the same as asking for
$s_{ij}=s_{g(i)g(j)}$ for any $i,j\in X$, so theorem 2.1 applies.

\bigskip

In general, commutation of $v$ with $d$ doesn't really depend on the precise value of coefficients $d_{ij}$. What matters is whether various pairs of coefficients $d_{ij}$ and $d_{kl}$ are equal or not. The same happens for an arbitrary coaction $v$ and an arbitrary matrix $d$.

\begin{theo}
Let $v$ be a coaction on a finite set $X$ and let $d\in M_X(\c )$ be a matrix. Consider the decomposition $d=\sum c\, d_c$, where for $c\in\c$ the matrix $d_c$ is defined by $(d_c)_{ij}=1$ if $d_{ij}=c$ and $(d_c)_{ij}=0$ if not, and where the sum is over nonzero terms. Then $(v_{ij})$ commutes with $d$ if and only if it commutes with all matrices $d_c$.
\end{theo}

\begin{proof}
We follow the proof in \cite{sms}. The magic biunitarity condition shows that the multiplication $M:\delta_i\otimes \delta_j\mapsto \delta_i\delta_j$ and comultiplication $C:\delta_i\mapsto \delta_i\otimes \delta_i$ intertwine $v^{\otimes 2}$ and $v$. Their iterations $M^{(k)}$ and $C^{(k)}$ intertwine $v^{\otimes k}$ and $v$, so the following operator commutes with $v$.
$$d^{(k)}=M^{(k)}d^{\otimes k}C^{(k)}=\sum_{c\in\c}c^k d_c$$

Let $S$ be the set of complex numbers $c$ such that $d_c\neq 0$. Consider the function $f:S\to\c$ given by $c\to c$ for any $c$. This function separates points of $S$, and by the Stone-Weierstrass theorem the subalgebra of $\c (S)$ generated by $f$ must be $\c (S)$ itself. In particular for any $e\in S$ the Dirac mass at $e$ is a linear combination of powers of $f$.
$$\delta_e=\sum_{k}\lambda_kf^k=\sum_{k}\lambda_k \left(\sum_{c\in S}c^k\delta_c\right)
=\sum_{c\in S}\left(\sum_{k}\lambda_kc^k\right)\delta_c$$

The corresponding linear combination of matrices $d^{(k)}$ is given by the following formula.
$$\sum_k\lambda_kd^{(k)}=\sum_k\lambda_k \left(\sum_{c\in S}c^kd_c\right)
=\sum_{c\in S}\left(\sum_{k}\lambda_kc^k\right)d_c$$

Dirac masses being linearly independent, in the first formula all coefficients in the right term are 0, except for the coefficient of $\delta_e$, which is 1. Thus the right term in the second formula is $d_e$. It follows that $d_e$ is is the algebra $End(v)$ of operators commuting with $v$.
\end{proof}

Another useful decomposition of $d$ is the spectral decomposition, in case $d$ is self-adjoint. If so, the $\c^*$-algebra generated by $d$ is spanned by spectral projections, and commutation with $v$ gives invariant subspaces, by using the following simple fact.

\begin{theo}
Let $v$ be a coaction and let $K$ be a linear subspace of $\c (X)$. The matrix $(v_{ij})$ commutes with the projection onto $K$ if and only if $v(K)\subset K\otimes H$.
\end{theo}

\begin{proof}
Let $P$ be the projection onto $K$. For any point $k\in X$ we have the following formula.
$$v(P(\delta_k))=v\left(\sum_i P_{ik}\delta_i\right) 
=\sum_{ij}P_{ik}\delta_j\otimes v_{ji}
=\sum_j\delta_j\otimes \sum_iv_{ji}P_{ik}$$

On the other hand the linear map $(P\otimes id)v$ is given by a similar formula.
$$(P\otimes id)(v(\delta_k))=\sum_i P(\delta_i)\otimes v_{ik}
=\sum_{ij}\delta_j\otimes P_{ji}v_{ik}=\sum_j\delta_j\otimes \sum_iP_{ji}v_{ik}$$

It follows that $vP=(P\otimes id)v$ is equivalent to the following conditions, for any $j,k$.
$$\sum_iv_{ji}P_{ik} =\sum_iP_{ji}v_{ik}$$

In other words, the equality of linear maps $vP=(P\otimes id)v$ is equivalent to the equality of products of square matrices $vP=Pv$, and the conclusion follows.
\end{proof}

The universal Hopf $\c^*$-algebra coacting on $X$ is constructed by Wang in \cite{wa2}. Its quotient by the relations $vd=dv$ is a universal object for the notion of coaction we are interested in.

\begin{theo}
Let $X$ be a finite set and let $d\in M_X(\c )$ be a complex matrix. Consider the universal $\c^*$-algebra $H(X,d)$ defined with generators $v_{ij}$ with $i,j\in X$ and with the relations making $v=(v_{ij})$ a magic biunitary matrix commuting with $d$. Then $H(X,d)$ is a Hopf $\c^*$-algebra and $v$ is a coaction of it on $X$.
\end{theo}

\begin{proof}
The universal $\c^*$-algebra exists because its generators $v_{ij}$ are projections, whose norms are bounded by $1$. Call it $H$. Consider the following matrix with coefficients in $H\otimes H$.
$$w_{ij}=\sum_kv_{ik}\otimes v_{kj}$$

Since $v$ is a magic biunitary, $w$ is magic biunitary as well. We multiply to the right by $d$.
$$(wd)_{ij}=\sum_kv_{ik}\otimes (vd)_{kj}=\sum_kv_{ik}\otimes (dv)_{kj}
=\sum_{kl}d_{kl}v_{ik}\otimes v_{lj}$$

On the other hand, the product $dw$ is given by the same formula.
$$(dw)_{ij}=\sum_l(dv)_{il}\otimes v_{lj}=\sum_l (vd)_{il}\otimes v_{lj}
=\sum_{kl}d_{kl}v_{ik}\otimes v_{lj}$$

Thus $w$ is a magic biunitary commuting with $d$, and the formula $\Delta(v_{ij})=\sum v_{ik}\otimes v_{kj}$ defines a $\c^*$-morphism. Also, the identity matrix is a magic biunitary commuting with $d$, so the formula  $\varepsilon (v_{ij})=\delta_{ij}$ defines a $\c^*$-morphism. Consider now the transpose matrix $v^t$, whose coefficients are viewed as elements of the opposite algebra $H^{op}$. Then $v^t$ is a magic biunitary, and commutation with $d$ follows from the following computation in $H$.
$$v^td=v^*d=v^*d(vv^*)=v^*(dv)v^*=v^*(vd)v^*=(v^*v)dv^*=dv^*=dv^t$$

Thus $S(v_{ij})=v_{ji}$ defines a $\c^*$-morphism, and all conditions in theorem 1.2 are satisfied.
\end{proof}

In order to cut off unwanted complexity, the very first condition to be put on $(X,d)$ is quantum homogeneity. This is a condition which goes under various names -- in \cite{sms} it is called quantum transitivity -- stating that the algebra of functions fixed by $v$ reduces to $\c$. This is the same as saying that the associated planar algebra or subfactor is irreducible.

\begin{defi}
We say that $(X,d)$ is quantum homogeneous if $v(f)=f\otimes id$ implies that $f$ is a constant function, where $v$ is the universal coaction of $H(X,d)$.
\end{defi}

The main example is when $(X,d)$ is homogeneous, meaning that its symmetry group $G$ acts transitively. (That is, for any $i,j\in X$ there is a permutation $\sigma :X\to X$ which preserves $d$, such that $\sigma (i)=j$.) Indeed, transitivity of $G$ is equivalent to the fact that $w(f)=f\otimes id$ implies that $f$ is constant, where $w$ is the corresponding coaction of $\c (G)$, and by using the universal property of $v$ we get that $(X,d)$ is quantum homogeneous.

We don't know if the converse holds, namely if quantum homogeneous implies homogeneous, but we have the following useful criterion here.

\begin{theo}
Let $d\in M_X(\c )$ be a complex matrix. For any complex number $c$ the characteristic function of $\{ i\mid d_{ii}=c\}$ is fixed by the universal coaction of $H(X,d)$.

In particular if $(X,d)$ is quantum homogeneous all diagonal entries of $d$ must be equal.
\end{theo}

\begin{proof}
Let $f$ be the characteristic function. By using theorem 2.2 we may assume that $d$ is a 0--1 matrix and that $c=1$. Consider the decomposition of the 0--1--(--1) matrix $d-1$.
$$d-1=(d-1)_1-(d-1)_{-1}$$

Theorem 2.2 shows that $v$ commutes with $P=(d-1)_{-1}$. Since --1 values can appear only on the diagonal of $d-1$, the 0--1 matrix $P$ is diagonal and $1-P$ is the projection onto $\c f$. Thus $v$ commutes with the projection onto $\c f$ and theorem 2.3 applies.
$$v(\c f)\subset (\c f)\otimes H(X,d)$$

We can write $v(f)=f\otimes a$. From the natural condition we get $\Sigma (f)a=\Sigma (f)1$, so $a=1$.
\end{proof}

\section{Colored semi-oriented graphs}

We are interested in spaces $(X,d)$ which are quantum homogeneous. Theorem 2.5 shows that all diagonal entries of $d$ must be equal. By substracting a scalar multiple of the identity we may assume that $d$ is zero on the diagonal.

This framework is still too general, and we don't have further results at this level. For the rest of the paper we make the quite natural assumption that $d$ is self-adjoint.

So, assume that $d$ is self-adjoint and has 0 on the diagonal. We call vertices the elements of $X$. For any $i$ and $j$ consider the complex number $c=d_{ij}=\overline{d}_{ji}$. If $c=0$ we do nothing, if $c$ is real we draw the edge $ij$ and color it $c$, if the imaginary part of $c$ is positive we draw the oriented edge $ij$ and color it $c$, and if the imaginary part of $c$ is negative we draw the oriented edge $ji$ and color it $\overline{c}$. We get a picture, that we call colored semi-oriented graph.

\begin{defi}
A colored semi-oriented graph $X$ is a finite graph with all edges colored and possibly oriented, such that an oriented edge and a non-oriented one cannot have same color. The choice of colors is not part of $X$.
\end{defi}

For each color $c$ consider the semi-oriented graph $X_c$ obtained by removing all edges having color different from $c$, then by considering that remaining edges are no longer colored $c$. These are graphs and oriented graphs, called color components of $X$.

The incidency matrices of a graph and of an oriented graph are by definition the matrices $d$ and $s$ in examples 2.2 and 2.3. The incidency matrices of colored components of $X$ are called incidency matrices of $X$ and are denoted $d_c$. Those corresponding to graphs are 0--1 matrices, those corresponding to oriented graphs are 0--$i$--$(-i)$ matrices. They are all self-adjoint. The collection of all incidency matrices determines $X$.

\begin{defi}
Associated to $X$ is the universal Hopf $\c^*$-algebra $H(X)$ coacting on the set of vertices, such that the matrix of coefficients $(v_{ij})$ commutes with all incidency matrices of $X$.
\end{defi}

Both existence and uniqueness follow from theorems 2.2 and 2.4. Note that for graphs our notion of coaction is different from Bichon's notion \cite{bic}, where $H$ has a coaction on the algebra of functions on the set of edges, compatible with $v$ in some natural sense. In the $H=\c (G)$ case both our notion of coaction and Bichon's coincide with the usual notion for groups. In the general case they are different. For instance when $X$ is the complete graph with 4 vertices definition 3.2 gives Wang's algebra in \cite{wa2}, which is infinite dimensional, while Bichon's universal construction produces the algebra $\c (S_4)$, cf. comments before proposition 3.3 in \cite{bic}.

Let $X$ and $Y$ be colored semi-oriented graphs having the same vertex set. By analogy with usual symmetry groups, we say that $H(X)$ is ``bigger'' than $H(Y)$ if there exists a Hopf $\c^*$-algebra morphism $H(X)\to H(Y)$ mapping coefficients of the universal coaction on $X$ to corresponding coefficients of the universal coaction on $Y$. This happens precisely when for any coaction $v$ of a Hopf $\c^*$-algebra on the vertex set we have that commutation of $(v_{ij})$ with the incidency matrices of $Y$ implies commutation  of $(v_{ij})$ with the incidency matrices of $X$.

It follows from definitions that bigger and smaller imply equal. Note also that the algebra of functions on the symmetry group $G(X)$ is smaller than $H(X)$.

\begin{theo}
Let $X$ be a colored semi-oriented graph.

(a) Removing a color component increases $H(X)$.

(b) Reversing orientation in a color component doesn't change $H(X)$.

(c) Forgetting orientation in a color component increases $H(X)$.

(d) Identifying two different colors, assumed to color same type of edges, increases $H(X)$.

(e) Adding all missing edges, unoriented and colored with a new color, doesn't change $H(X)$.
\end{theo}

\begin{proof}
Removing a color component means removing a commutation relation in definition of $H$, and we get (a). At level of incidency matrices reversing orientation is given by $d\to -d$, and we get (b). It is enough to prove (c) for an oriented graph. Here the incidency matrix has 0--1 decomposition of type $s=is_i-is_{-i}$, and the new incidency matrix is given by $d=s_i+s_{-i}$. If $v$ commutes with $s$ it must commute with both matrices $s_{\pm i}$, so it commutes with $d$ as well. It is enough to prove (d) for graphs, oriented or not. In both cases commutation of $v$ with two matrices is replaced by commutation with their sum or difference, and we get (d).

From (a) we get that adding all missing edges decreases $H(X)$. Thus in (e) it is enough to prove that $H(X)$ increases. For, let $X_b$ be obtained from $X$ by forgetting all orientations and identifying all colors, say with a black color. The incidency matrices of $X_b$ and of the new color component, say a white color component, are related by the following formula.
$$d_{b}+d_{w}=\begin{pmatrix}
0&1&\ldots&1\cr
1&0&\ldots&1\cr
\ldots&\ldots&\ldots&\ldots\cr
1&1&\ldots&0\end{pmatrix}$$

Let $v$ be the universal coaction of $H(X)$. If $\i$ denotes the matrix filled with 1 then left and right multiplication by $\i$ is making sums on rows and columns, and the magic biunitary condition shows that both $v\i$ and $\i v$ are equal to $\i$. Thus $v$ commutes with the matrix on the right. On the other hand we know from (c) and (d) that $H(X_b)$ is bigger than $H(X)$, so $v$ commutes with $d_b$. It follows that $v$ commutes with $d_w$ and we are done.
\end{proof}

\begin{theo}
If $X$ is a graph and $X^c$ is its complement then $H(X)=H(X^c)$.
\end{theo}

\begin{proof}
From (e) we get that both $H(X)$ and $H(X^c)$ are equal to $H$ of the complete graph having color components $X$ and $X^c$.
\end{proof}

It is convenient to give names to pictures. Most geometric objects are metric spaces, and in the non-oriented case we can use (e) plus the following consequence of theorem 2.2.

\begin{theo}
To any finite metric space we associate the colored graph $X$ obtained by drawing edges between all pairs of points and coloring them with their lenghts. Then $H(X)$ is isomorphic to the universal Hopf $\c^*$-algebra coacting on the space in a co-isometric way, in the sense that matrix of coefficients $(v_{ij})$ commutes with the distance matrix.
\end{theo}

In general, when drawing a picture of a geometric object what happens is that the symmetry group of the picture, regarded as a graph, is equal to the symmetry group of the rigid object, regarded as a metric space. The same is true for associated Hopf algebras.

\begin{exam}
Simplex vs simplex.
\end{exam}

Consider the $n$-simplex, viewed as metric space with $n$ points. Consider also the usual picture of the $n$-simplex, viewed as a graph with $n$ vertices and $(^n_2)$ edges. The distance matrix of the metric space is proportional to the incidency matrix of the graph.

Thus $H$ of the $n$-simplex metric space is equal to $H$ of the $n$-simplex graph.

\begin{exam}
Cube vs cube.
\end{exam}

Consider the cube, viewed as metric space with 8 points. Consider also the usual picture of a cube, viewed as a graph with 8 vertices and 12 edges. The 0--1 decomposition of the distance matrix of the metric space is as follows, where $a$ is the lenght of the side.
$$d=a\, d_{a}+\sqrt{2}a\, d_{\sqrt{2}a}+\sqrt{3}a\, d_{\sqrt{3}a}$$

The 0--1 matrix $d_a$ is the incidency matrix of the graph. Its square counts 2-loops on the graph, and we get $d_a^2=3+2d_{\sqrt{2}a}$. Thus a coaction $v$ on the set of 8 vertices commuting with $d_1$ must commutes with $d_{\sqrt{2}a}$, and from theorem 3.1 we get that $v$ commutes with $d_{\sqrt{3}a}$ as well. Thus $v$ must commute with $d$. Conversely, if $v$ commutes with $d$ then theorem 2.2 shows that $v$ commutes with $d_a$. Commutation with $d$ is equivalent to commutation with $d_a$.

Thus $H$ of the cube metric space is equal to $H$ of the cube graph.

\begin{exam}
Polygon vs polygon.
\end{exam}

Consider the regular $n$-gon, viewed as a metric space. Consider also the usual picture of the regular $n$-gon, viewed as a graph. The 0--1 decomposition of the distance matrix of the metric space is as follows, where $a$ is the lenght of the sides and $b<c<\ldots$ are the lenghts of various diagonals.
$$d=a\, d_a+b\, d_b+c\, d_c+\ldots$$

The 0--1 matrix $d_a$ is the incidency matrix of the graph. Counting 2-loops gives $d_a^2=2+d_b$, so commutation with $d_a$ is equivalent to commutation with both $d_a$ and $d_b$. The picture shows that $d_ad_b=d_a+d_c$, $d_ad_c=d_b+d_d$ and so on, and by induction we get that commutation with $d_a$ is equivalent to commutation with all 0--1 components of $d$.

Thus $H$ of the regular $n$-gon metric space is equal to $H$ of the regular $n$-gon graph.

\begin{theo}
Let $X$ be a graph which is quantum homogeneous, in the sense that $v(f)=f\otimes 1$ implies that $f$ is constant, where $v$ is the universal coaction commuting with the incidency matrix. For any $l\geq 2$ the number of $l$-loops based at a vertex is independent of the vertex.
\end{theo}

\begin{proof}
This follows from the fact that the diagonal entry $d^l(p,p)$ of the $l$-th power of the incidency matrix counts $l$-loops at $p$, and in theorem 2.5 one can replace $d$ by any of its powers.
\end{proof}

\section{Cyclic and dihedral groups}

The simplest example of oriented graph is the oriented $n$-gon. This graph has vertices $1,\ldots ,n$ and an oriented edge joins $i$ and $i+1$ for any $i$, with $i$ taken modulo $n$.

\begin{theo}
The Hopf $\c^*$-algebra associated to the oriented $n$-gon is the algebra of functions on the cyclic group $\z_n$. The Poincar\'e series is given by
$$f(z)=1+\frac{z}{1-nz}$$
and its convergence radius is $1/n$.
\end{theo}

\begin{proof}
The oriented $n$-gon $X$ has the following real incidency matrix.
$$r=\begin{pmatrix}
0&1&0&\dots &0\cr
0&0&1&\dots &0\cr
\dots&\dots&\dots&\dots &\dots\cr
0&0&0&\dots &1\cr
1&0&0&\dots &0\cr
\end{pmatrix}$$

This is a permutation matrix. If $v$ is the coaction of $H(X)$, commutation of $v$ with $r$ says that $v$ must be of the following special form.
$$v=\begin{pmatrix}
v_1&v_2&v_3&\dots &v_n\cr
v_n&v_1&v_2&\dots &v_{n-1}\cr
\dots&\dots&\dots&\dots &\dots\cr
v_3&v_4&v_5&\dots &v_2\cr
v_2&v_3&v_4&\dots &v_1\cr
\end{pmatrix}$$

Since $v$ is a magic biunitary, the elements $v_i$ form a partition of the unity of $H(X)$. In particular they commute, so $H(X)$ is commutative. Thus $H(X)$ is the algebra of functions on the usual symmetry group of $X$, which is the cyclic group $\z_n$.

The Poincar\'e series is computed by using theorem 1.4. The unit of $\z_n$ has $n$ fixed points, and the other $n-1$ elements, none.
$$f(z)=\frac{1}{n}\left( n-1+\frac{1}{1-nz}\right) 
=1+\frac{1}{n}\left( -1+\frac{1}{1-nz}\right) =1+\frac{1}{n}\cdot\frac{nz}{1-nz}$$

This is equal to the function in the statement.
\end{proof}

In the non-oriented case, the Hopf $\c^*$-algebra associated to the $n$-gon is infinite dimensional if $n=4$, and equal to the algebra of functions on the dihedral group $D_n$ if $n\neq 4$. For $n\leq 4$ this is known from Wang's paper \cite{wa2} and for $n\geq 5$ this is proved in \cite{sms}.

We extend now the $n\neq 4$ result to a larger class of cyclic graphs.

\begin{defi}
We say that a graph $X$ is cyclic if its automorphism group contains a copy of the cyclic group $\z_n$, where $n$ is the number of vertices.
\end{defi}

If $X$ is cyclic, one can choose a vertex and label it $0$, then label the other vertices $1,\ldots ,n-1$ such that $\z_n$ acts by $g(k)=g+k$, with both $g$ and $k$ modulo $n$.

It is useful to keep in mind the following interpretation. Vertices of $X$ are $n$-th roots of unity in the complex plane, counted counterclockwise starting with $1$, and edges are segments joining vertices. The graph is cyclic if the $2\pi /n$ rotation of the plane leaves invariant the picture.

\begin{defi}
If $X$ is a cyclic graph with vertices labeled $0,1,\ldots ,n-1$ we define numbers $e(k)$ by $e(k)=1$ if $0$ and $k$ are connected by an edge and $e(k)=0$ if not. Then
$$Q(z)=e(1)z+e(2)z^2+\ldots +e(n-2)z^{n-2}+e(n-1)z^{n-1}$$
is a polynomial which doesn't depend on the choice of the vertex $0$.
\end{defi}

This is a sum of monomials which is symmetric with respect to $z^{n/2}$, in the sense that $Q$ has degree at most $n-1$ and the coefficient of $z^k$ is equal to the coefficient of $z^{n-k}$, for any $k$. Any sum of monomials which is symmetric with respect to $z^{n/2}$ is of this form.

The simplest cyclic graph is the $n$-gon, corresponding to $Q(z)=z+z^{n-1}$.

\begin{theo}
Let $X$ be a cyclic graph with $n\neq 4$ vertices and consider the associated polynomial $Q$. Let $w$ be a primitive $n$-th root of unity. If the numbers
$$Q(1),\, Q(w),\, Q(w^2 ),\ldots ,Q(w^{[n/2]})$$
are distinct then $H(X)$ is the algebra of functions on $D_n$. The Poincar\'e series is
$$f(z)=1+\frac{z}{2}\left(\frac{1}{1-\varepsilon z}+\frac{1}{1-nz}\right)$$
where $\varepsilon =1$ if $n$ is odd and $\varepsilon =2$ if $n$ is even. The convergence radius is $1/n$.
\end{theo}

This is an extension of the result for the $n$-gon. Indeed, for the $n$-gon having vertices at roots of unity $1,w,w^2,\ldots ,w^{n-1}$ the number $Q(w^k)=w^k+w^{n-k}$ for $w^k$ above the $x$-axis is twice the projection of $w^k$ on the $x$-axis, which decreases when $k$ increases.

\begin{proof}
This follows proofs in \cite{wa2}, \cite{sms}. We use the Vandermonde formula.
$$\begin{pmatrix}
0     &e(1)  &e(2)  &\ldots&e(n-1)\cr
e(n-1)&0     &e(1)  &\ldots&e(n-2)\cr
\dots&\dots&\dots&\dots&\dots\cr
e(1)  &e(2)  &e(3)  &\ldots&0
\end{pmatrix}
\begin{pmatrix}1\cr w^k\cr w^{2k}\cr\dots\cr w^{(n-1)k}\end{pmatrix}
=Q(w^k)\begin{pmatrix}1\cr w^k\cr w^{2k}\cr\dots\cr w^{(n-1)k}\end{pmatrix}$$

The eigenvalues $Q(w^k)$ being distinct, we get the list of invariant subspaces of $d$.
$$\c 1,\,\c\xi\oplus\c\xi^{n-1},\, \c\xi^2\oplus\c\xi^{n-2},\ldots$$

Here $\xi =(w^i)$ and the last subspace has dimension 1 or 2, depending on the parity of $n$. We claim that $v$ is given by formulae $F_k$ of following type, with $k=1,2,\ldots ,n-1$.
$$v(\xi^k)=\xi^k\otimes a^k+\xi^{n-k}\otimes b^k$$

This follows by applying many times theorem 2.3. For $n=2$ we take $b=0$, we define $a$ by $F_1$, and $F_k$ follows by induction. For $n\geq 3$ we define $a,b$ by $F_1$. Taking the square and cube of $F_1$ gives $ab=-ba$ and $ab^2=ba^2=0$. With these relations, $F_k$ follows by induction.

Applying $*$ to $F_1$ and comparing with $F_{n-1}$ gives $a^*=a^{n-1}$ and $b^*=b^{n-1}$. Together with $ab^2=0$ this gives $abb^*a^*=0$. Thus $ab=ba=0$, and in particular $H(X)$ is commutative. On the other hand $H(X)$ depends only on spectral projections of $d$, so it must be the same for all graphs in the statement. With the $n$-gon we get $H(X)=\c (D_n)$.

For $n=1$ the Poincar\'e series is computed by using theorem 1.4.
$$f(z)=\frac{1}{1-z}=1+z\cdot\frac{1}{1-z}=
1+\frac{z}{2}\left(\frac{1}{1-z}+\frac{1}{1-z}\right)$$

For $n=2$ the group $G$ has two elements. The identity is in $G_2$ and the other element is in $G_0$. This gives the formula in the statement.
$$f(z)=\frac{1}{2}\left( 1+\frac{1}{1-2z}\right) =\frac{1-z}{1-2z}
=1+\frac{z}{1-2z}=1+\frac{z}{2}\left(\frac{1}{1-2z}+\frac{1}{1-2z}\right)$$

For $n\geq 3$ odd the group $G$ has $2n$ elements. The identity is in $G_n$, the $n-1$ rotations are in $G_0$ and the $n$ symmetries are in $G_1$. Theorem 1.4 applies and gives the result.
\begin{eqnarray*}
f(z)&=&\frac{1}{2n}\left( n-1+\frac{n}{1-z}+\frac{1}{1-nz}\right)\\
&=&1+\frac{1}{2n}\left(\left(\frac{n}{1-z}-n\right)+\left(\frac{1}{1-nz}-1\right)\right)\\
&=&1+\frac{1}{2n}\left( \frac{nz}{1-z}+\frac{nz}{1-nz}\right)
\end{eqnarray*}

For $n\geq 4$ even the group $G$ has $2n$ elements as well. The identity is in $G_n$ and the $n-1$ rotations are in $G_0$. There are $n$ more elements, namely the symmetries, half of them being in $G_0$ and half of them being in $G_2$. We apply theorem 1.4.
\begin{eqnarray*}
f(z)&=&\frac{1}{2n}\left( 3n/2-1+\frac{n/2}{1-2z}+\frac{1}{1-nz}\right)\\
&=&1+\frac{1}{2n}\left(\left(\frac{n/2}{1-2z}-n/2\right)+\left(\frac{1}{1-nz}-1\right)\right)\\
&=&1+\frac{1}{2n}\left( \frac{nz}{1-2z}+\frac{nz}{1-nz}\right)
\end{eqnarray*}

The proof of theorem 4.2 is now complete.
\end{proof}

We apply now theorem 4.2 to graphs with small number of vertices. We call first and second 9-star the graphs corresponding to the following polynomials, with $e=1,2$.
$$Q(z)=z+z^{1+e}+z^{8-e}+z^8$$

The 8-spoke wheel is the graph corresponding to $Q(z)=z+z^4+z^7$.

\begin{coro}
If $X$ is a graph with $n\in \{1,\ldots ,9\} -\{ 4\}$ vertices having symmetry group $D_n$ then $H(X)=\c (D_n)$. The list of such graphs is as follows.

(a) $n$-gons with $n\neq 4$ and their complements.

(b) $8$-spoke wheel and its complement.

(c) $9$-stars.
\end{coro}

\begin{proof}
Let $X$ be as in the statement. By relabeling 048372615 vertices of the first 9-star we see that its complement is the second 9-star. Thus the list of graphs is closed under complementation. It is enough to show that $X$ or $X^c$ appears in the list. By replacing $X\to X^c$ we may assume that the valence $k$ of vertices is smaller than $(n-1)/2$. Thus $n\geq 2k+1$.

For $k=0$ the possible graphs are the point, the 2 points and the 3 points. These are the 1-gon, the complement of the 2-gon, and the complement of the 3-gon, all 3 in the list.

For $k=1$ we have $n\geq 3$. On the other hand the graph must be a union of segments, so the only solution is the 2 segments. But here $n=4$.

For $k=2$ we have the $n$-gons with $n\geq 5$, all of them in the list.

For $k=3$ we have $n\geq 7$. The cases $n=7,9$ are excluded, because the incidency matrix must have $3n/2$ values of 1 above the diagonal. In the $n=8$ case the graph corresponds to a polynomial of form $Q(z)=z^a+z^4+z^{8-a}$ with $a=1,2,3$. For $a=1$ this is the 8-spoke wheel, in the list, for $a=2$ we get the 2 tetrahedra, not dihedral, and for $a=3$ we can relabel vertices 03614725 and we get the 8-spoke wheel again.

For $k=4$ we have $n=9$. Consider the associated $Q$ polynomial.
$$Q(z)=z^a+z^b+z^{9-b}+z^{9-a}$$

Here $a\neq b$ are from $\{ 1,2,3,4\}$. Since $a,b$ are not both equal to 3, one of them, say $a$, is prime with 3, and by relabeling vertices $0,a,2a,\ldots ,7a$ we can assume $a=1$. For $b=2,3$ we get the 9-stars and for $b=4$ we can relabel vertices 048372615 and we get the first 9-star.

Theorem 4.2 applies to $n$-gons with $n\neq 4$. For the 8-spoke wheel we have $Q(w^k)=w^k+(-1)^k+w^{-k}$ and computation gives the distinct numbers $3$, $\pm 1$ and $-1\pm\sqrt{2}$, so theorem 4.2 applies as well. Let $X$ be a 9-star and assume that we have an equality of the form $Q(w^k)=Q(w^l)$ with $0\leq k<l\leq 4$. The numbers $Q(w^n)/2$ are sums of numbers $\cos(2n\pi/9)$, and we see that the only equalities between two such sums are those of form $x+x=x+x$. In particular we must have $\cos(2k\pi/9)=\cos(2l\pi/9)$. But this is impossible because $\cos(2x\pi/9)$ is decreasing on $[0,4]$. Thus theorem 4.2 applies to both 9-stars and we are done.
\end{proof}

\section{Tannaka-Galois duality}

An arbitrary Hopf $\c^*$-algebra with a faithful coaction $v:\c (X)\to\c (X)\otimes H$ is quite an abstract object. However, a useful description is obtained after classifying its irreducible corepresentations, together with their fusion rules. In this section we present a categorical and topological approach to this problem, by using Woronowicz's Tannakian duality \cite{wo2} and the spin planar algebra $P(X)$ constructed by Jones in \cite{j1}, \cite{j2}. The main result will be a Tannaka-Galois type correspondence between pairs $(H,v)$ and subalgebras $P\subset P(X)$, somehow in the spirit of the correspondence found by Kodiyalam, Landau and Sunder in \cite{kls}.

Recall from section 1 that the main problem is to decompose tensor powers of $v$. Consider the Hilbert space where the $m$-th tensor power of $v$ acts.
$$P_m(X)=\c (X)^{\otimes m}$$

In planar calculus both the input and the output are written in a $2\times m$ matrix form, and the first thing to be done is to write elements of $P_m(X)$ in such a way.

\begin{defi}
Each $m$-fold tensor product of Dirac masses at points of $X$ is identified with a $2\times m$  matrix with coefficients in $X$, in the following way.
$$\delta_{i_1}\otimes \ldots\otimes\delta_{i_m}=\begin{pmatrix}i_1 & i_1 &i_2&i_2&i_3&\ldots\cr i_m & i_m &i_{m-1}&\ldots&\ldots&\ldots \end{pmatrix}$$

That is, we take the sequence of $m$ points, we duplicate each entry, then we put it in the $2\times m$ matrix, clockwise starting from top left.
\end{defi}

We recall now some basic notions from Jones' planar algebra formalism \cite{j1}.

A box is a rectangle in the plane. It is convenient to assume that sides of the box are parallel to the real axis and imaginary axis. We say that a box $X$ is at left of a box $Y$ if the horizontal sides of $X$ and $Y$ are on the same lines, and if $X$ is at left of $Y$. Same for on top.

A $m$-box is a box with $2m$ marked points, $m$ of them on the lower side and $m$ of them on the upper side. If $x\in P_m(X)$ is a tensor product of Dirac masses, written in loop form, we can put indices on marked points in the obvious way. We say that $x$ is in the box.

Let $m_1,\ldots ,m_k$ and $n$ be positive integers. Let $T$ be a picture consisting of an output $n$-box, containing an input $m_i$-box for each $i$, together with some non-crossing strings outside the input boxes. These are strings connecting pairs of marked points, plus a finite number of closed strings, called circles. Strings connecting input and output points are assumed to connect odd-numbered points to odd-numbered points and even to even, when numbering  points on each $k$-box $1,2,\ldots ,2k$ clockwise starting from top left. Such a picture, or rather its planar isotopy class, is called $(m_1,\ldots ,m_k,n)$-tangle, or just $m$-tangle.

Planar tangles act on tensors in the following way. See Jones \cite{j1}, \cite{j2}.

\begin{defi}
Each $(m_1,\ldots ,m_k,n)$-tangle $T$ defines a multilinear map
$$P_{m_1}(X)\otimes\ldots\otimes P_{m_k}(X)\to P_n(X)$$
in the following way. If $x_1,\ldots ,x_k$ and $y$ are $m_1,\ldots ,m_k$-fold and $n$-fold tensor products of Dirac masses written in loop form, put each $x_i$ in the $m_i$-input box of $T$ and $y$ in the output $n$-box of $T$. Strings of $T$ join now indices, and the number $(x_1,\ldots ,x_k,y)^T$ is defined to be $1$ if all strings join pairs of equal indices and $0$ if not. Define
$$T(x_1\otimes\ldots\otimes x_k)=\beta^{c(T)}\sum_y(x_1,\ldots ,x_k,y)^Ty$$
where the sum is over all $n$-fold tensor products of Dirac masses, $\beta$ is the number of elements of $X$ and $c(T)$ is the number of closed circles of $T$.
\end{defi}

The planar calculus for tensors is quite simple and doesn't really require diagrams. It suffices to imagine that the way various indices appear, travel around and dissapear is by following some obvious strings connecting them. Some illustrating examples.

\begin{exam}
Identity, multiplication, inclusion.
\end{exam}

The identity $1_m$ is the $(m,m)$-tangle having vertical strings only. The solutions of $(x,y)^{1_m}=1$ are pairs of the form $(x,x)$, so $1_m$ acts by the identity.
$$1_m\begin{pmatrix}j_1 & \ldots & j_m\cr i_1 & \ldots & i_m\end{pmatrix}=\begin{pmatrix}j_1 & \ldots & j_m\cr i_1 & \ldots & i_m\end{pmatrix}$$

The multiplication $M_m$ is the $(m,m,m)$-tangle having 2 input boxes, one on top of the other, and vertical strings only. It acts in the following way.
$$M_m\left( 
\begin{pmatrix}j_1 & \ldots & j_m\cr i_1 & \ldots & i_m\end{pmatrix}
\otimes\begin{pmatrix}l_1 & \ldots & l_m\cr k_1 & \ldots & k_m\end{pmatrix}
\right)=
\delta_{j_1k_1}\ldots \delta_{j_mk_m}
\begin{pmatrix}l_1 & \ldots & l_m\cr i_1 & \ldots & i_m\end{pmatrix}$$

The formula $xy=M_m(x\otimes y)$ defines an associative multiplication of $P_m(X)$.

The inclusion $I_m$ is the $(m,m+1)$-tangle which looks like $1_m$, but has one more vertical string, at right of the input box. Given $x$ written in loop form, solutions of $(x,y)^{I_m}=1$ are elements $y$ obtained from $x$ by adding to the right a vector of the form $(^l_l)$.
$${I_m}\begin{pmatrix}j_1 & \ldots & j_{m}\cr i_1 & \ldots & i_{m}\end{pmatrix}=
\sum_l\begin{pmatrix}j_1 & \ldots & j_{m} & l\cr i_1 & \ldots & i_{m}& l\end{pmatrix}$$

This shows that $I_m$ is an inclusion of algebras, and that various $I_m$ are compatible with each other. The inductive limit of the algebras $P_m(X)$ is a graded algebra, denoted $P(X)$.

\begin{exam}
Expectation, Jones projection.
\end{exam}

The expectation $U_m$ is the $(m+1,m)$-tangle which looks like $1_m$, but has one more string, connecting the extra 2 input points, that we suppose to be both at right of the input box.
$$U_m
\begin{pmatrix}j_1 & \ldots &j_m& j_{m+1}\cr i_1 & \ldots &i_m& i_{m+1}\end{pmatrix}=
\delta_{i_{m+1}j_{m+1}}
\begin{pmatrix}j_1 & \ldots & j_{m}\cr i_1 & \ldots & i_{m}\end{pmatrix}$$

This shows that $U_m$ is a bimodule morphism with respect to $I_m$.

The Jones projection $E_m$ is a $(0,m+2)$-tangle, having no input box. There are $m$ vertical strings joining the first $m$ upper points to the first $m$ lower points, counting from left to right. The remaining upper 2 points are connected by a semicircle, and the remaining lower 2 points are also connected by a semicircle. We can apply $E_m$ to the unit of $\c$.
$$E_m(1)=\sum\begin{pmatrix}i_1 & \ldots &i_m& j&j\cr i_1 & \ldots &i_m&k&k\end{pmatrix}$$

The elements $e_m=(\# X)^{-1}E_m(1)$ are projections, and define a representation of the infinite Temperley-Lieb algebra of index $\# X$ on the inductive limit algebra $P(X)$.

\begin{exam}
Rotation.
\end{exam}

The rotation $R_m$ is the $(m,m)$-tangle which looks like $1_m$, but the first 2 input points are connected to the last 2 output points, and the same happens at right.
$$R_5=\begin{matrix}
\hskip 0.3mm\Cap \ |\ |\ |\ |\hskip -0.5mm |\cr
|\hskip -0.5mm |\hskip 10.3mm |\hskip -0.5mm |\cr
\hskip -0.3mm|\hskip -0.5mm |\ |\ |\ |\ \hskip -0.1mm\Cup
\end{matrix}$$

The action of $R_m$ is best described in terms of Dirac masses.
$$R_m\left( \delta_{i_1}\otimes \ldots\otimes\delta_{i_m}\right)
=\delta_{i_2}\otimes \ldots\otimes\delta_{i_m}\otimes\delta_{i_1}$$

Thus $R_m$ acts by an order $m$ linear automorphism of $P_m(X)$, also called rotation.

\bigskip

Multiplications, inclusions, expectations, Jones projections, rotations generate in fact the set of all tangles, with the gluing operation described below.

\begin{defi}
Let $T$ be a $(m_1,\ldots ,m_k,n)$-tangle and let $S$ be a $m_i$-tangle with $m_i$ among $m_1,\ldots ,m_k$. The composition $TS$ is obtained by superposing the output box of $S$ on the input $m_i$-box of $T$, after isotoping so that marked points match, then by removing the common boundary. The colored planar operad $\p$ is the set of all tangles, with this gluing operation.
\end{defi}

The composition $U_mI_m$ consists of $1_m$ plus a floating circle, and by using the above formulae for actions of $1_m$, $I_m$ and $U_m$ we get the following equality.
$$(U_mI_m)(x)=U_m(I_m(x))=(\# X)\, x=1_m^\circ (x)$$

In general, composition of tangles corresponds to composition of maps. We have a morphism from $\p$ to the colored operad of multilinear maps between spaces $P(X)$, called action of $\p$ on $P(X)$. This action commutes in some natural sense with the involution of $P_m(X)$.
$$\begin{pmatrix}j_1 & \ldots & j_m\cr i_1 & \ldots & i_m\end{pmatrix}^*
=\begin{pmatrix}i_1 & \ldots & i_m\cr j_1 & \ldots & j_m\end{pmatrix}$$

This means that $P(X)$ is a $\c^*$-planar algebra, called spin planar algebra. See Jones \cite{j1}.

\begin{defi}
The graded linear space $P(X)$ together with the action of $\p$ and with the involution $*$ is called spin planar algebra associated to $X$.
\end{defi}

Let $v:\c (X)\to\c (X)\otimes H$ be a coaction. Consider the $m$-th tensor power of $v$.
$$v^{\otimes m}:P_m(X)\to P_m(X)\otimes H$$

Computation using the magic biunitarity condition shows that each $v^{\otimes m}$ is a $\c^*$-morphism. See \cite{pl} for details. Let $P_m$ be the fixed point algebra of $v^{\otimes m}$.
$$P_m=\{ x\in P_m(X)\mid v^{\otimes m}(x)=x\otimes 1\}$$

Consider the rotation $R_m$. Rotating, then applying $v^{\otimes m}$, then rotating backwards by $R_m^{-1}$ is the same as applying $v^{\otimes m}$, then rotating each $m$-fold product of coefficients of $v$. Thus the elements obtained by rotating, then applying $v^{\otimes m}$, or by applying $v^{\otimes m}$, then rotating, differ by a sum of tensor products of Dirac masses tensor commutators in $H$.
$$v^{\otimes m}R_m(x)-(R_m\otimes id)v^{\otimes m}(x)\in P_m(x)\otimes [H,H]$$

Let $h$ be the Haar functional and consider the conditional expectation $\phi_m =(id\otimes h)v^{\otimes m}$ onto the fixed point algebra $P_m$. The square of the antipode being the identity, $h$ is a trace, so it vanishes on commutators. Thus $R_m$ commutes with $\phi_m$.
$$\phi_m R_m=R_m\phi_m$$

The commutation relation $\phi_nT=T\phi_m$ holds in fact for any $(m,n)$-tangle $T$. These tangles are called annular, and proof in \cite{pl} is by verification on generators of the annular category. In particular we get $\phi_nT\phi_m=T\phi_m$ for any $T$, so the annular category is contained in the suboperad $\p^\prime\subset \p$ consisting of tangles $T$ satisfying the following condition, where $\phi =(\phi_m)$ and $i(.)$ is the number of input boxes.
$$\phi T\phi^{\otimes i(T)}=T\phi^{\otimes i(T)}$$

On the other hand multiplicativity of $v^{\otimes m}$ gives $M_m\in\p^\prime$. Since $\p$ is generated by multiplications and annular tangles, it follows that $\p^\prime =P$. Thus for any tangle $T$ the corresponding multilinear map between spaces $P_m(X)$ restricts to a multilinear map between spaces $P_m$. In other words, the action of $\p$ restricts to $P$ and makes it a subalgebra of $P(X)$.

\begin{defi}
The sequence of spaces of fixed points of $v^{\otimes m}$, together with the restriction of the action of $\p$ is called $\c^*$-planar algebra associated to $v$.
\end{defi}

Consider pairs $(H,v)$ where $v:\c (X)\to \c(X)\otimes H$ is a faithful coaction. A morphism $(H,v)\to (K,w)$ is a Hopf $\c^*$-algebra morphism $H\to K$ sending $v_{ij}\to w_{ij}$. Isomorphism means morphisms in both senses, and this is the notion used in the statement below.

\begin{theo}
If $Q$ is a $\c^*$-planar subalgebra of $P(X)$ there is a unique pair $(H,v)$ with $v:\c (X)\to \c(X)\otimes H$ faithful coaction whose associated $\c^*$-planar algebra is $Q$.
\end{theo}

\begin{proof}
This will follow by applying Woronowicz's Tannakian duality \cite{wo2} to the annular category over $Q$. This is constructed as follows. Let $n,m$ be positive integers. To any element $T_{n+m}\in Q_{n+m}$ we associate a linear map $L_{nm}(T_{n+m}):P_n(X)\to P_m(X)$ in the following way.
$$L_{nm}\left(\begin{matrix}|\ |\ |\cr T_{n+m}\cr |\ |\ |\end{matrix}\right):
\left(\begin{matrix}|\cr a_n\cr |\end{matrix}\right)
\to \left(\begin{matrix}
\hskip 1.5mm |\hskip 3.0mm |\hskip 3.0mm \cap\cr
\ \ T_{n+m}\hskip 0.0mm  |\cr
\hskip 1.9mm |\hskip 1.2mm |\hskip 3.2mm |\hskip2.2mm |\cr
a_n|\hskip 3.2mm |\hskip 2.2mm |\cr
\hskip 2.1mm\cup \hskip3.5mm |\hskip 2.2mm |
\end{matrix}\right)$$

That is, we consider the planar $(n,n+m,m)$-tangle having an small input $n$-box, a big input $n+m$-box and an output $m$-box, with strings as on the picture of the right. This defines a certain multilinear map $P_n(X)\otimes P_{n+m}(X)\to P_m(X)$. Now we put $T_{n+m}$ in the big input box. What we get is a linear map $P_n(X)\to P_m(X)$. This is called $L_{nm}$.

The above picture corresponds to $n=1$ and $m=2$. This is illustrating whenever $n\leq m$, suffices to imagine that in the general case all strings are multiple.

If $n>m$ there are $n+m$ strings of $a_n$ which connect to the $n+m$ lower strings of $T_{n+m}$, and the remaining $n-m$ ones go to the upper right side and connect to the $n-m$ strings on top right of $T_{n+m}$. Here is the picture for $n=2$ and $m=1$.
$$L_{nm}\left(\begin{matrix}|\ |\ |\cr T_{n+m}\cr |\ |\ |\end{matrix}\right):
\left(\begin{matrix}|\ |\cr a_n\cr |\ |\end{matrix}\right)
\to \left(\begin{matrix}
\hskip 3mm |\hskip 3.7mm \Cap\cr
T_{n+m}|\hskip -0.5mm |\cr
\hskip 4.7mm ||\hskip 1.8mm |\hskip -0.5mm |\hskip -0.5mm |\cr
\hskip 4.7mm a_n|\hskip -0.5mm |\hskip -0.5mm |\cr
\hskip 7.5mm \Cup \hskip -0.2mm|\cr
\hskip 10.1mm|
\end{matrix}\right)$$

This problem with two cases $n\leq m$ and $n>m$ can be avoided by using an uniform approach, with discs with marked points instead of boxes. See Jones \cite{j3}.

Consider the linear spaces formed by such maps.
$$Q_{nm}=\{ L_{nm}(T_{n+m}):P_n(X)\to P_m(X)\mid T_{n+m}\in Q_{n+m}\}$$

Pictures show that these spaces form a tensor $\c^*$-subcategory of the tensor $\c^*$-category of linear maps between tensor powers of the Hilbert space $H=\c (X)$. If $j$ is the antilinear map from $\c (X)$ to itself given by $j(^i_i)=(^i_i)$ for any $i$, then the elements $t_j (1)$ and $t_{j^{-1}} (1)$ constructed by Woronowicz in \cite{wo2} are both equal to the unit of $Q_2=Q_{02}$. In other words, the tensor $\c^*$-category has conjugation, and Tannakian duality in \cite{wo2} applies.

We get a pair $(H,v)$ consisting of a unital Hopf $\c^*$-algebra $H$ and a unitary corepresentation $v$ of $H$ on $\c (X)$, such that the following equalities hold, for any $m,n$.
$$Hom(v^{\otimes m},v^{\otimes n})=Q_{mn}$$

We prove that $v$ is a magic biunitary. We have $Hom(1,v^{\otimes 2})=Q_{02}=Q_2$, so the unit of $Q_2$ must be a fixed vector of $v^{\otimes 2}$. But $v^{\otimes 2}$ acts on the unit of $Q_2$ in the following way.
$$v^{\otimes 2}(1)
=v^{\otimes 2}\left( \sum_i \begin{pmatrix}i&i\cr i&i\end{pmatrix}\right) 
=\sum_{ikl}\begin{pmatrix}k&k\cr l&l\end{pmatrix}\otimes v_{ki}v_{li}
=\sum_{kl}\begin{pmatrix}k&k\cr l&l\end{pmatrix}\otimes (vv^t)_{kl}$$

From $v^{\otimes 2}(1)=1\otimes 1$ ve get that $vv^t$ is the identity matrix. Together with the unitarity of $v$, this gives the following formulae.
$$v^t=v^*=v^{-1}$$

Consider the Jones projection $E_1\in Q_3$. After isotoping $L_{21}(E_1)$ looks as follows.
$$L_{21}\left( \Bigl| \begin{matrix}\cup\cr\cap\end{matrix}\right) :
\begin{pmatrix} \,|\ |\cr {\ }^i_j{\ }^i_j\cr \,|\ |\end{pmatrix}\,\to\,
\begin{pmatrix}\hskip -5.8mm |\cr {\ }^i_j{\ }^i_j\supset\cr \hskip -5.8mm |\end{pmatrix}
=\,\delta_{ij}\begin{pmatrix}\,|\cr {\ }^i_i\cr \,|\end{pmatrix}$$

In other words, the linear map $M=L_{21}(E_1)$ is the multiplication $\delta_i\otimes\delta_j\to\delta_{ij}\delta_i$.
$$M\begin{pmatrix}i&i\cr j&j\end{pmatrix}=\delta_{ij}
\begin{pmatrix}i\cr i\end{pmatrix}$$

We have $M\in Q_{21}=Hom(v^{\otimes 2},v)$, so the following elements of $\c (X)\otimes H$ are equal.
$$(M\otimes id)v^{\otimes 2}\left(\begin{pmatrix}i&i\cr j&j\end{pmatrix}\otimes 1\right)
=(M\otimes id)\left(\sum_{kl}\begin{pmatrix}k&k\cr l&l\end{pmatrix}
\otimes v_{ki}v_{lj}\right)
=\sum_{k}\begin{pmatrix}k\cr k\end{pmatrix}\delta_k\otimes v_{ki}v_{kj}$$
$$v(M\otimes id)\left(\begin{pmatrix}i&i\cr j&j\end{pmatrix}\otimes 1\right)
=v\left(\delta_{ij}\begin{pmatrix}i\cr i\end{pmatrix}\delta_i\otimes 1\right)
=\sum_k\begin{pmatrix}k\cr k\end{pmatrix}\delta_k\otimes\delta_{ij}v_{ki}$$

Thus $v_{ki}v_{kj}=\delta_{ij}v_{ki}$ for any $i,j,k$. With $i=j$ we get $v_{ki}^2=v_{ki}$, and together with the formula $v^t=v^*$ this shows that all entries of $v$ are self-adjoint projections. With $i\neq j$ we get $v_{ki}v_{kj}=0$, so projections on each row of $v$ are orthogonal to each other. Together with $v^t=v^{-1}$ this shows that each row of $v$ is a partition of unity with self-adjoint projections.

The antipode is given by the formula $(id\otimes S)v=v^*$. But $v^*$ is the transpose of $v$, so we can apply $S$ to the formulae saying that rows of $v$ are partitions of unity, and we get that columns of $v$ are also partitions of unity. Thus $v$ is a magic biunitary.

Consider the planar algebra $P$ associated to $v$. We have the following equalities.
$$Hom(1,v^{\otimes n})=P_n$$

Thus $P_n=Q_n$ for any $n$ and this proves the existence assertion.

As for uniqueness, let $(K,w)$ be another pair corresponding to $Q$. The functorial properties of Tannakian duality give a morphism $f:(H,v)\to (K,w)$. Since morphisms increase spaces of fixed points we have the following inclusions.
$$Q_k=Hom(1,v^{\otimes k})\subset Hom(1,w^{\otimes k})=Q_k$$

We must have equality for any $k$, and by using Frobenius reciprocity and a basis of coefficients of irreducible corepresentations we see that $f$ must be an isomorphism on this basis and we are done. This is a standard argument, see for instance lemma 5.3 in \cite{subf}.
\end{proof}

\begin{coro}
If $X$ is a graph with $n$ vertices and no edges the planar algebra associated to $H(X)$ is the Temperley-Lieb algebra $TL(n)$. For $n\geq 4$ the Poincar\'e series is
$$f(z)=\frac{1-\sqrt{1-4z}}{2z}$$
with convergence radius $1/4$. If $n=1,2,3$ we have $H(X)=\c (D_n)$ and the Poincar\'e series is
$$f(z)=1+\frac{z}{2}\left(\frac{1}{1-\varepsilon z}+\frac{1}{1-nz}\right)$$
where $\varepsilon =1$ if $n=1,3$ and $\varepsilon =2$ if $n=2$, with convergence radius $1/n$.
\end{coro}

\begin{proof}
First assertion follows from theorem 5.1, because the universal object $H(X)$ must correspond to the universal object $TL(n)\subset P(X)$. The formula for the Poincar\'e series of the Temperley-Lieb algebra is well-known. As for the last assertion, this follows either from theorem 4.2 or from the well-known fact that $TL(n)$ corresponds to $D_n$ for $n=1,2,3$.
\end{proof}

\section{Exchange relations and multi-simplexes}

Let $X$ be a finite set and $d\in M_X(\c )$ be a self-adjoint matrix with indices in $X$. This matrix can be viewed as a 2-box in the spin planar algebra $P(X)$, by using the canonical identification between $M_X(\c )$ and the algebra of 2-boxes $P_2(X)$.
$$(d_{ij})\leftrightarrow \sum_{ij} d_{ij}\begin{pmatrix}i&i\cr j&j\end{pmatrix}$$

The planar algebra generated by $d$ is the smallest linear subspace of $P(X)$ containing $d$ and stable by the action of the colored planar operad $\p$.

\begin{theo}
The planar algebra associated to $H(X,d)$ is equal to the planar algebra generated by $d$, viewed as a $2$-box in the spin planar algebra $P(X)$.
\end{theo}

\begin{proof}
Let $P$ be the planar algebra associated to $H(X,d)$ and let $Q$ be the planar algebra generated by $d$. The action of $v^{\otimes 2}$ on $d$ viewed as a 2-box is given by the following formula.
$$v^{\otimes 2}\left(\sum_{ij} d_{ij}\begin{pmatrix}i&i\cr j&j\end{pmatrix}\right)
=\sum_{ijkl} d_{ij}\begin{pmatrix}k&k\cr l&l\end{pmatrix}\otimes v_{ki}v_{lj}
=\sum_{kl}\begin{pmatrix}k&k\cr l&l\end{pmatrix}\otimes (vdv^t)_{kl}$$

Since $v$ is a magic biunitary commuting with $d$ we have $vdv^t=dvv^t=d$. This means that $d$, viewed as a 2-box, is in the algebra $P_2$ of fixed points of $v^{\otimes 2}$. Thus $Q\subset P$.

For $P\subset Q$ we apply theorem 5.1. Let $(K,w)$ be the pair whose associated planar algebra is $Q$. The same computation with $w$ at the place of $v$ shows that $w$ commutes with $d$. Thus the universal property of $H(X,d)$ gives a Hopf $\c^*$-algebra morphism $H(X,d)\to K$ sending $v_{ij}\to w_{ij}$. Since morphisms increase spaces of fixed points we have the following inclusions.
$$P_k=Hom(1,v^{\otimes k})\subset Hom(1,w^{\otimes k})=Q_{0k}=Q_k$$

It follows that $P\subset Q$ and we are done.
\end{proof}

The planar algebra generated by $d$ is also constructed by Curtin in \cite{cu}.

With $d=0$ this gives an alternative proof of corollary 5.1, stating that a simplex corresponds to a Temperley-Lieb algebra. This is the same as a Fuss-Catalan algebra on 1 color.

The fact that a product of 2 simplexes gives a Fuss-Catalan algebra on 2 colors is known from \cite{sms}. This is a consequence of a straightforward computation in \cite{fc}, using the presentation result of Bisch-Jones \cite{bj1}. We give here a simpler proof, which works for $s$ colors. This uses theorem 6.1 and the presentation result of Landau \cite{la}.

The Fuss-Catalan algebra on $s$ colors is presented by a sequence of $2$-boxes $p_1,\ldots ,p_s$ satisfying relations in theorem 3 in \cite{la}. In the particular case of a subalgebra of the spin planar algebra $P(X)$ we can identify 2-boxes with matrices.
$$\begin{matrix}|\, |\cr (p_i)\cr |\, |\end{matrix}=
\sum_{ab}p^i_{ab}
\begin{matrix}|\ |\cr (^a_b{\,}^a_b)\cr |\ |\end{matrix}$$

Landau's relations for $FC(n_1,\ldots ,n_s)$ are as follows. First is a standard one.

$$\sum_{ab}p^i_{ab}
\begin{matrix}|\ |\cr (^a_b{\,}^a_b) \cr |\ |\end{matrix}
=\sum_{ab}\overline{p^i_{ab}}
\begin{matrix}|\ |\cr(^b_a{\,}^b_a)\cr |\ |\end{matrix}
=\sum_{ab}p^i_{ab}
\begin{matrix}|\ |\cr(^b_a{\,}^b_a)\cr |\ |\end{matrix}$$

In the second relation we have two formulae, one of them involving indices.

$$\sum_{ab}p^i_{ab}
\begin{matrix}|\ |\cr(^a_b{\,}^a_b)\cr\cup\end{matrix}=
\begin{matrix}\bigl|\, \bigl|\cr\cup\end{matrix}
\hskip 2cm
\sum_{ab}p^i_{ab}
\begin{matrix}\ |\ \ \cap\cr (^a_b{\,}^a_b) \bigl|\cr \ |\ \ \cup\end{matrix}
=n_1\ldots n_{i-1}\biggl|$$

Third is the exchange relation. This must hold for any $i\geq j$.
$$\sum_{ab\alpha\beta}p_{ab}^ip_{\alpha\beta}^j
\begin{matrix}
\, |\ |\ \ \ \hskip0.3mm |\cr
(^a_b{\,}^a_b) \ \ |\ \cr
|\ \backslash \ \ /\cr
| (^\alpha_\alpha{\,}^\beta_\beta )\cr
\hskip 0.3mm |\ /\ \ \backslash
\end{matrix}=
\sum_{ab\alpha\beta}p_{ab}^ip_{\alpha\beta}^j
\begin{matrix}
\hskip 0.3mm |\ \backslash\ \ /\cr
| (^\alpha_\alpha{\,}^\beta_\beta )\cr
|\ / \ \ \backslash\cr
(^a_b{\,}^a_b) \ \ |\ \cr
\, |\ |\ \ \ \hskip0.3mm |
\end{matrix}$$

The missing index $n_s$ is implicit from the embedding into $P(X)$. Observe that we use a reverse labeling of projections $p_i$. This is for the comment after definition 6.1 to make sense.

These relations can be reformulated by using the spin algebra specific conventions.
$$p^i_{ab}=\overline{p^i_{ba}}=p^i_{ba}\hskip 1cm \sum_bp_{ab}^i=1\hskip 1cm\sum_ap_{aa}^i=n_1\ldots n_{i-1}\hskip 1cm p_{ab}^ip_{b\beta}^j=p_{ab}^ip_{a\beta}^j$$

Here all relations must hold for all choices of indices, except for the last one, which is subject to the condition $i\geq j$. The solution is given by the following notion.

\begin{defi}
A $(n_1,n_2,\ldots ,n_s)$-simplex is a colored graph $X$ obtained in the following way. The vertex set is the product of sets $Y_i$ having $n_i$ elements each
$$X=Y_1\times\ldots \times Y_s$$
and the edges of $X$ are all possible edges, colored in the following way.

(1) For $z_1\neq t_1$ the edge between any $(z_1,\ldots )$ and any $(t_1,\ldots )$ is colored $\varepsilon_1$.

(2) For $z_2\neq t_2$ the edge between any $(z_1,z_2,\ldots )$ and any $(z_1,t_2,\ldots )$ is colored $\varepsilon_2$.

\dots

($s$) For $z_s\neq t_s$ the edge between any $(z_1,\ldots ,z_{s-1},z_s)$ and any $(z_1,\ldots ,z_{s-1},t_s)$ is colored $\varepsilon_s$.
\end{defi}

It is useful to keep in mind the following ``metric'' interpretation. We have finite sets $Y_1,\ldots ,Y_s$ and a decreasing sequence of $s$ infinitesimals.
$$\varepsilon_1>>\varepsilon_2>>\ldots >>\varepsilon_s$$

(1) Start with the simplex $Y_1$ with distance $\varepsilon_1$.

(2) Put a copy of the simplex $Y_2$ with distance $\varepsilon_2$ at each vertex of $Y_1$.

\dots

($s$) Put a copy of $Y_s$ with distance $\varepsilon_s$ at each vertex of $Y_{s-1}$.

This is a kind of metric space, by using the infinitesimal summing conventions $\varepsilon_i+\varepsilon_j=\varepsilon_i$ for $i>j$. It is also possible to start with a galaxy having solar systems and so on.

\begin{theo}
If $X$ is a $(n_1,n_2,\ldots ,n_s)$-simplex the planar algebra associated to $H(X)$ is the Fuss-Catalan algebra with $s$ colors and indices $n_1,n_2,\ldots ,n_s$. For generic indices $n_i\geq 4$ the Poincar\'e series is given by the Bisch-Jones formula
$$f(z)=\sum_{k=0}^{\infty}\frac{1}{sk+1}\begin{pmatrix}(s+1)k\cr k\end{pmatrix}z^k$$
and its convergence radius is $s^s/(s+1)^{(s+1)}$.
\end{theo}

\begin{proof}
As already mentioned, this is just a Hopf $\c^*$-algebra reformulation of a particular case of results of Bisch and Jones \cite{bj1} and Landau \cite{la}. We use the following identification.
$$M_X(\c )=M_{Y_1}(\c )\otimes\ldots \otimes M_{Y_s}(\c )$$

For any $i$ let $1_i$ be $n_i\times n_i$ matrix having 1 on the diagonal and 0 outside, and let $\i_i$ be the  $n_i\times n_i$ matrix having 1 everywhere. The incidency matrices of $X$ are as follows.
\begin{eqnarray*}
d_1&=&(\i_1-1_1)\otimes \i_2\otimes\i_3\otimes\ldots \otimes \i_s\\
d_2&=&1_1\otimes (\i_2-1_2)\otimes \i_3\otimes\ldots \otimes \i_s\\
&\dots &\ \\
d_s&=&1_1\otimes\ldots\otimes 1_{s-1}\otimes (\i_s-1_s)
\end{eqnarray*}

We have $d_i=e_{i}-e_{i+1}$ for any $i$, where $e_{s+1}=1$ and where $e_i$ are the following operators.
\begin{eqnarray*}
e_1&=&\i_1\otimes \i_2\otimes\ldots \otimes \i_s\\
e_2&=&1_1\otimes \i_2\otimes\ldots \otimes \i_s\\
&\dots &\ \\
e_{s}&=&1_1\otimes\ldots\otimes 1_{s-1}\otimes \i_s
\end{eqnarray*}

Commutation with all matrices $d_i$ is equivalent to commutation with all matrices $e_i$. Thus the planar algebra associated to $H(X)$ is the planar algebra generated by matrices $e_i$, viewed as 2-boxes in the spin planar algebra $P(X)$. On the other hand each $e_i$ is a scalar multiple of a certain projection $p_i$, and we can replace matrices $e_i$ by projections $p_i$.
$$p_i=\frac{e_i}{n_i\ldots n_s}$$

The coefficients of $p_i$ are given by the following formula.
$$p^i_{(a_1\ldots a_s)(b_1\ldots b_s)}=
\frac{\delta_{a_1b_1}\ldots\delta_{a_{i-1}b_{i-1}}}{n_i\ldots n_s}$$

We verify now Landau's relations. First one is clear from the following identity.
$$\delta_{a_1b_1}\ldots\delta_{a_{i-1}b_{i-1}}
=\overline{\delta_{b_1a_1}\ldots\delta_{b_{i-1}a_{i-1}}}
=\delta_{b_1a_1}\ldots\delta_{b_{i-1}a_{i-1}}$$

The two formulae forming the second relation are clear as well.
$$\sum_{b_1\ldots b_s}
\frac{\delta_{a_1b_1}\ldots\delta_{a_{i-1}b_{i-1}}}{n_i\ldots n_s}=1\hskip 1cm
\sum_{a_1\ldots a_s}
\frac{\delta_{a_1a_1}\ldots\delta_{a_{i-1}a_{i-1}}}{n_i\ldots n_s}
=n_1\ldots n_{i-1}$$

For any $i\geq j$ the exchange relation follows from the following formula.
$$\delta_{a_1b_1}\ldots\delta_{a_{i-1}b_{i-1}}
\delta_{b_1\beta_1}\ldots\delta_{b_{j-1}\beta_{j-1}}
=\delta_{a_1b_1}\ldots\delta_{a_{i-1}b_{i-1}}
\delta_{a_1\beta_1}\ldots\delta_{a_{j-1}\beta_{j-1}}$$

By \cite{la} the planar algebra generated by projections $p_i$ is isomorphic to the Fuss-Catalan algebra on $s$ colors. Surjectivity follows from positivity of the Markov trace, see \cite{fc}. The formula of $f$ is from Bisch and Jones \cite{bj1}, and the convergence radius is well-known, and easy to compute by using the Stirling formula.
\end{proof}

For non-generic indices see Bisch and Jones \cite{bj1}, \cite{bj15}.

\begin{coro}
The planar algebra associated to two squares $\Box\Box$ is $FC(2,2,2)$.
\end{coro}

\begin{proof}
The square of the incidency matrix of $\Box\Box$ is $2$ times the identity plus the incidency matrix of the graph XX formed by missing diagonals. Theorem 3.1 shows that $H(\Box\Box)$ is equal to $H$ of the bicolored graph obtained by superposing a black copy of XX on a blue copy of $\Box\Box$, which is in turn is equal to $H$ of the $(2,2,2)$-simplex, which gives $FC(2,2,2)$.
\end{proof}

\section{Decomposable graphs}

If $Y$ and $Z$ are graphs we define a graph $Y\times Z$ in the following way. Vertices are pairs $(y,z)$ with $y$ vertex of $Y$ and $z$ vertex of $Z$ and there is an edge in between $(y_1,z_1)$ and $(y_2,z_2)$ if there are edges in $Y$ between $y_1$ and $y_2$ and in $Z$ between $z_1$ and $z_2$.

The simplest example of such a decomposition is $X_{2n}=X_n\times X_2$, where $X_m$ denotes the $m$-gon. In this identification $2k+e$ corresponds to $(k,e)$, where vertices of all $m$-gons are labeled counterclockwise $0,1,\ldots ,m-1$. On the other hand, the automorphism groups of these graphs are related by the formula $D_{2n}=D_n\times D_2$ if $n$ is odd. If $n$ is even this equality fails. This can be interpreted in the following way. For an $n$-gon $Y$ and a segment $Z$ the Hopf $\c^*$-algebras $H(Y\times Z)$ and $H(Y)\otimes H(Z)$ are isomorphic if and only if $n$ is odd.

If $Y$ and $Z$ are finite sets we use the following canonical identification.
$$\c (Y\times Z)=\c (Y)\otimes \c (Z)$$

We say that a graph is regular if all its vertices have the same valence.

\begin{theo}
Let $Y$ and $Z$ be finite connected regular graphs. If the sets of eigenvalues of their incidency matrices $\{\lambda\}$ and $\{\mu\}$ don't contain $0$ and satisfy the condition 
$$\left\{ \lambda_i/\lambda_j\right\}\cap
\left\{ \mu_k/\mu_l\right\} =\{ 1\}$$
then $H(Y\times Z)=H(Y)\otimes H(Z)$. The Poincar\'e series is the coefficientwise product of Poincar\'e series for $Y$ and for $Z$, and its convergence radius is the product of convergence radii.
\end{theo}

\begin{proof}
Let $\lambda_1$ be the valence of $Y$. Since $Y$ is connected $\lambda_1$ has multiplicity 1, and if $P_1$ is the orthogonal projection onto $\c 1$, the spectral decomposition of $d_Y$ is of the following form.
$$d_Y=\lambda_1\, P_1 +\sum_{i\neq 1}\lambda_i\, P_i$$

We have a similar formula for $d_Z$, where $Q_1$ is the projection onto $\c 1$.
$$d_Z=\mu_1\, Q_1 +\sum_{j\neq 1}\mu_j\, Q_j$$

It follows from definitions that the incidency matrix of $Y\times Z$ is $d_Y\otimes d_Z$.
$$d_{Y\times Z}=\sum_{i,j}\lambda_i\mu_j (P_i\otimes Q_j)$$

In this formula the projections form a partition of the unity and the scalars are distinct, so this is the spectral decomposition of $d_{Y\times Z}$. The universal coaction $V$ on $Y\times Z$ must commute with all spectral projections.
$$[V, P_i\otimes Q_j]=0$$

By summing over $i$ relations with $j=1$ we get that $V$ commutes with $1\otimes Q_1$. But this is the projection onto the algebra $\c (Y)\otimes \c 1$, so this algebra must be invariant. The corresponding restriction of $V$ is a coaction of $H(Y\times Z)$ on the spectrum of $\c (Y)\otimes \c 1$, and can be written by using a magic biunitary $w=(w_{xy})$ with indices in $Y$.
$$V\left( \delta_x\otimes 1 \right) =\sum_y \delta_y\otimes 1\otimes w_{yx}$$

The same argument gives a similar formula for the restriction to $\c 1\otimes \c (Z)$.
$$V\left( 1\otimes\delta_p \right) =\sum_q 1\otimes\delta_q\otimes u_{qp}$$

By multiplying these equalities in all possible ways we get two formulae for $V$.
$$V\left(\delta_x\otimes \delta_p\right) 
=\sum_{yq} \delta_y\otimes \delta_q\otimes w_{yx}u_{qp}
=\sum_{yq} \delta_y\otimes \delta_q\otimes u_{qp}w_{yx}$$

Thus $w$ and $u$ commute and we have $V=w\otimes u$, where $w,u,V$ are regarded now as corepresentations, and where $\otimes$ is the usual tensor product of corepresentations. Let $W$ and $U$ be the fundamental corepresentations of $H(Y)$ and $H(Z)$. The universal properties of $H(Y)$ and $H(Z)$ defines morphisms of Hopf $\c^*$-algebras $\phi$ and $\psi$ sending $W_{yx}\to w_{yx}$ and $U_{qp}\to u_{qp}$. The images of these maps commute, so by tensoring we get a map sending $W\otimes U\to V$.
$$\phi\otimes\psi :H(Y)\otimes H(Z)\to H(Y\times Z)$$

On the other hand the corepresentation $W\otimes U$ can be viewed as a magic biunitary with coefficients in $Y\times Z$, and the universal property of $H(Y\times Z)$ gives an inverse map for $\phi\otimes\psi$. The assertions on Poincar\'e series are clear.
\end{proof}

For graphs $Y\times Z$ with $n\leq 7$ vertices the only application is to the product of a triangle and a segment, which gives $D_6=D_3\times D_2$. For $n=8$ the only possible application is with a tetrahedron and a segment, where we use the conventions from section 3.

\begin{coro}
For the cube $X$ we have $H(X)=H(T)\otimes\c (\z_2)$, where $T$ is the graph with 4 vertices and no edges. The Poincar\'e series is given by
$$f(z)=1+\sum_{k=1}^{\infty}\frac{2^{k-1}}{k+1}\begin{pmatrix}2k\cr k\end{pmatrix}z^k$$
and has convergence radius $1/8$.
\end{coro}

\begin{proof}
By highlighting all diagonals of lenght $\sqrt{2}$ we get two graphs of tetrahedra, say $Y$ and $Y^\circ$. Let $Z$ be a segment, with vertices denoted $0$ and $1$. We have an identification $X=Y\times Z$, where vertices $p\in Y$ correspond to pairs $(p,0)$ and vertices $q\in Y^\circ$ correspond to pairs $(q^\circ,1)$, where $\circ$ is the symmetry with respect to the center of the cube. For $Z$ the eigenvalue set is $\{ -1,1\}$ and for $Y$ the eigenvalue set is $\{ -1,3\}$.
$$\begin{pmatrix}0&1&1&1\cr 1&0&1&1\cr 1&1&0&1\cr 1&1&1&0\end{pmatrix}
=3\left( \frac{1}{4}\begin{pmatrix}1&1&1&1\cr 1&1&1&1\cr 1&1&1&1\cr 1&1&1&1\end{pmatrix}\right)
-\left( 1-\frac{1}{4}\begin{pmatrix}1&1&1&1\cr 1&1&1&1\cr 1&1&1&1\cr 1&1&1&1\end{pmatrix}\right)$$

Theorem 7.1 applies. By using theorem 3.2 we can replace in the conclusion $Y$ with its complement $T$ and this gives the first assertion.

The Poincar\'e series $f$ is the coefficientwise product of the Poincar\'e series $g$ and $h$ of the graphs formed by 4 and by 2 points. Corollary 5.1 shows that $g$ is the Poincar\'e series of the Temperley-Lieb algebra of index 4, given by the following formula.
$$g(z)=\frac{1-\sqrt{1-4z}}{2z}=\sum_{k=0}^{\infty}\frac{1}{k+1}\begin{pmatrix}2k\cr k\end{pmatrix}z^k$$

The complement of 2 points is the 2-gon, so $h$ is the Poincar\'e series of the action of $\c (D_2)$ on the 2-gon, given by theorem 4.2. Multiplying coefficients gives the second assertion.
\end{proof}

\section{Graphs with small number of vertices}

The most interesting graphs are those which are quantum homogeneous, in the sense that the algebra of functions fixed by universal coaction is 1-dimensional. This is an a priori abstract notion, and we don't know yet how to check it on the picture. Homogeneous implies quantum homogeneous, and in fact, all examples we have are of this form.

Recall from section 4 that the 8-spoke wheel is a graph with 8 vertices and 12 edges which looks like a wheel with 8 spokes, with the convention that the middle point is not a vertex.

\begin{theo}
A graph with $n\leq 8$ is quantum homogeneous if and only if it is homogeneous. There are $38$ such graphs, and results in previous sections apply to all of them.

(a) point, $2$ points, segment, $3$ points, triangle.

(b) graphs with $n=4,5,6,7,8$ vertices with no edges, or with all possible edges.

(c) $2$, $3$ or $4$ segments, $2$ triangles, $2$ tetrahedra, plus their complements.

(d) $2$ squares and its complement.

(e) pentagon, heptagon and its complement.

(f) hexagon, octogon, $8$-spoke wheel, and their complements.

(g) cube and its complement.
\end{theo}

\begin{proof}
We will use many times theorem 3.4, to be called ``loop rule''. With $l=2$ this shows that a quantum homogeneous graph must be regular, so we restrict attention to regular graphs. Let $X$ be a regular graph with $n\leq 8$ vertices, say of valence $k$. It is enough to prove that $X$ or its complement $X^c$ is in the list. By replacing $X$ with $X^c$ we may assume that the valence of $X$ is smaller than the valence of $X^c$. Thus we may assume $n\geq 2k+1$.

For $k=0$ the graph is $n$ points, hence in (a) or (b).

For $k=1$ the only graphs are the 2, 3 or 4 segments, hence in (c).

For $k=2$ we have $n\geq 5$. The graph must be a union of $m$-gons with values of $m$ greater than 3. The loop rule with various $l$'s shows that the $m$'s must be the same for all $m$-gons. All such graphs are in (c), (d), (e) and (f).

For $k=3$ we have $n\geq 7$. The $n=7$ case is excluded, because the incidency matrix must have $3n/2$ values of 1 above the diagonal. For $n=8$ we apply the loop condition with $l=3$. Let $\Delta$ be the number of triangles based at each point, meaning half of the number of 3-loops. We have $\Delta\leq 3$, with equality if $X$ is 2 tetrahedra, which is in (c). On the other hand each triangle appears globally 3 times, so $8\Delta\in 3\z$. Thus the case left is $\Delta =0$.

Assume that 2 vertices have the same set of 3 neighbors. There are 3 edges emanating from these 2+3=5 points towards the remaining 8-5=3 ones. These remaining points having valence 3, they must form a triangle, contradiction.

Choose 1 point and consider its 3 neighbors. There are 6 edges emanating from these 3 neighbors towards the remaining 8-1-3=4 points. The above remark applied twice shows that 2 such edges must reach the same point and 2 other ones must reach some other point. There are 6-2-2=2 edges left, having to reach one of the 4-1-1=2 points left. There are two cases, depending on whether these 2 edges reach the same point or not. In both cases the rest of the graph is uniquely determined, so the 2 possible solutions must be the cube and the 8-spoke wheel. Pictures show that indeed they are. These are in (f) and (g), so we are done.
\end{proof}

We end up with a few remarks, most of them describing work in progress.

\begin{rema}
Quantum homogeneous graphs with $n\leq 8$ vertices.
\end{rema}

The results in previous sections apply to the above 38 graphs, and we
get the Poincar\'e series for each of them. It is possible to write
down a big table, with names and pictures of graphs, fusion rules,
principal graphs, and various descriptions of Poincar\'e series and
their coefficients. However, this table is not very enlightening,
because of too many small index accidents.

As a conclusion here, it is more convenient to state a classification
result.

Consider the planar algebra associated to an arbitrary graph
$(X,d)$. This planar algebra can be
defined alternatively as (1) the one corresponding to
$H(X,d)$ via Tannaka-Galois duality, or (2) the subalgebra generated by
the incidency matrix $d$, viewed as a 2-box in the spin
planar algebra $P(X)$. See theorem 6.1.

We can say that a graph $(X,d)$ is a Fuss-Catalan graph if the planar
algebra associated to it is a Fuss-Catalan algebra. This
Fuss-Catalan algebra depends on the number $s$ of colors, and on the
$s$ values of the indices, that we can omit or not. For instance when $s=1$ we
say that $(X,d)$ is a Temperley-Lieb graph, or that it is a $TL(n)$
graph, where $n$ is the index.

We also say that $(X,d)$ is dihedral is $H(X,d)=\c (D_n)$, where
$n$ is the number of vertices of $X$. The value $n$ can be
specified or not, for instance $D_n$-graph means $H(X,d)=\c (D_n)$.

With these definitions, the above 38 graphs fall into three
classes. First is the Fuss-Catalan series of graphs. Second is the dihedral series, with $D_4$ missing, and with $D_n$
overlapping with $TL(n)$ for $n=1,2,3$. Third are the cube and its
complement, corresponding to a ``tensor product'' between $TL(2)$ and
$TL(4)$. This follows from results in previous sections.

\begin{rema}
Quantum homogeneous graphs with $n=9$ vertices.
\end{rema}

A graph with $9$ vertices is quantum homogeneous if and only if it is homogeneous. The list of such graphs is as follows.

(a) Fuss-Catalan graphs, namely $9$ points, $9$-simplex, $3$ triangles and its complement.

(b) Dihedral graphs, namely  $9$-gon and its complement, $9$-stars.

(c) Discrete torus. This is the product of a triangle with itself.

This enumeration can be obtained by applying theorem 3.4 with
$l=2,3,4$, but the proof is quite long, and
won't be given here.

We don't know how to compute the Poincar\'e series of the discrete
torus. Note that theorem 7.1 doesn't apply, because the 2
eigenvalue sets in the statement are equal.

\begin{rema}
Colored or oriented graphs with small number of vertices.
\end{rema}

One can try to classify graphs with small number of vertices which are
not colored / colored and not oriented / oriented / semi-oriented,
in a way similar to the one in remark 8.1 (corresponding to the
not colored and not oriented case). For instance in the colored and
not oriented case such a classification is obtained in \cite{sms}, for
graphs with $n\leq 7$ vertices.

Let us just say here that using results in this paper (1) in the colored and not oriented case it
is possible to compute Poincar\'e series up to $n=8$, but there is one
graph left : the 2 rectangles, and (2) in the not colored and oriented
case it is possible to compute Poincar\'e series up
to $n=6$, but there is one graph left : the 2 oriented triangles.

That is, so far the classification work stumbles into 3 graphs: the
discrete torus, the 2 rectangles and the 2 oriented triangles. We
believe that the discrete torus can be studied by using a vertex model, and that
the $X\to 2X$ operation can be investigated by using free
product techniques, but we don't have yet results in this sense.

\begin{rema} 
Loop condition and homogeneity.
\end{rema}

The ``loop condition'' in theorem 3.4 doesn't imply that the graph is
homogeneous, as explained to us by E. Ghys. Let $P$ be a projective
plane, meaning that $P$ consists of points and lines such that each
two lines cross and each two points are on a line. The incidency graph
$X(P)$ has points and lines as vertices, and edges are drawn only
between points and lines, when the point is on the line. The axioms of
$P$ show that $X(P)$ satisfies the loop condition, but if $P$ is one
of the non-desarguian projective planes then $X(P)$ is not
homogeneous.

We don't know if these graphs are quantum homogeneous or not.

\end{document}